\documentclass[11pt,a4paper]{article}

\usepackage[]{hyperref}
\usepackage{amsmath, amsthm, amssymb,bm,dcolumn}
\usepackage{graphicx}
\usepackage[square,sort,comma,numbers]{natbib}

\theoremstyle{plain}

\title{Comment on ``Hilbert's Sixth Problem: Derivation of Fluid Equations via Boltzmann's Kinetic Theory'' by Deng, Hani, and Ma}

\author{Shan Gao
\\Research Center for Philosophy of Science and Technology, 
\\ Shanxi University, Taiyuan 030006, P. R. China
\\ E-mail:  \href{mailto:gaoshan2017@sxu.edu.cn}{gaoshan2017@sxu.edu.cn}.}

\begin{document}

\maketitle

\begin{abstract}
Deng, Hani, and Ma [\href{https://arxiv.org/abs/2503.01800}{arXiv:2503.01800}] claim to resolve Hilbert’s Sixth Problem by deriving the Navier-Stokes-Fourier equations from Newtonian mechanics via an iterated limit: a Boltzmann-Grad limit (\(\varepsilon \to 0\), \(N \varepsilon^{d-1} = \alpha\) fixed) yielding the Boltzmann equation, followed by a hydrodynamic limit (\(\alpha \to \infty\)) to obtain fluid dynamics. Though mathematically rigorous, their approach harbors two critical physical flaws. First, the vanishing volume fraction (\(N \varepsilon^d \to 0\)) confines the system to a dilute gas, incapable of embodying dense fluid properties even as \(\alpha\) scales, rendering the resulting equations a rescaled gas model rather than a true continuum. Second, the Boltzmann equation’s reliance on molecular chaos collapses in fluid-like regimes, where recollisions and correlations invalidate its derivation from Newtonian dynamics. These inconsistencies expose a disconnect between the formalism and the physical essence of fluids, failing to capture emergent, density-driven phenomena central to Hilbert’s vision. We contend that the Sixth Problem remains open, urging a rethink of classical kinetic theory’s limits and the exploration of alternative frameworks to unify microscale mechanics with macroscale fluid behavior.
\end{abstract}

\section{Introduction}

Hilbert’s Sixth Problem, posed in 1900 as part of his seminal list of mathematical challenges, calls for the axiomatization of physics, particularly the rigorous derivation of continuum equations—such as the Navier-Stokes-Fourier system—from the laws of microscopic particle dynamics governed by Newtonian mechanics \cite{Hilbert}. This problem lies at the intersection of mathematics, physics, and philosophy, demanding not only formal proofs but also a coherent physical narrative connecting atomic-scale interactions to macroscopic fluid behavior. Over a century later, despite advances in kinetic theory and statistical mechanics, a complete resolution remains elusive, with debates persisting over the validity and scope of proposed derivations.

In their recent work, Deng, Hani, and Ma \cite{DengHaniMa} claim to address this challenge by deriving the Navier-Stokes-Fourier equations from Newtonian mechanics through a two-step process involving Boltzmann’s kinetic theory. Their approach, framed on the periodic torus \(\mathbb{T}^d\) (\(d = 2\) or \(3\)), proceeds as follows: 
\begin{itemize}
    \item \textbf{Step 1 (Newton to Boltzmann):} Starting with \(N\) hard-sphere particles of diameter \(\varepsilon\) on \(\mathbb{T}^d\), they employ the Boltzmann-Grad limit (\(N \to \infty\), \(\varepsilon \to 0\), with \(N \varepsilon^{d-1} = \alpha\) fixed) to derive the Boltzmann equation. This step enforces molecular chaos—a statistical independence assumption—to suppress recollisions and correlations, ensuring the one-particle density \(n(t, x, v)\) satisfies the irreversible Boltzmann equation.
    \item \textbf{Step 2 (Boltzmann to Navier-Stokes):} By scaling the collision frequency \(\alpha \to \infty\) (hydrodynamic limit), they drive \(n(t, x, v)\) to a local Maxwellian distribution, from which the Navier-Stokes-Fourier equations emerge via moment closure. The authors assert that this iterated limit (\(\varepsilon \to 0\), then \(\alpha \to \infty\)) bridges Newtonian mechanics to fluid dynamics, resolving Hilbert’s Sixth Problem \cite{DengHaniMa}.
\end{itemize}

While mathematically elegant, this derivation raises profound physical and conceptual questions. Central to Hilbert’s vision is the demand for a \textit{physically consistent} reduction—one that not only produces formal equations but also respects the ontology of fluids as dense, continuum media. Herein lies the crux of our critique: the authors’ framework, while valid for dilute gases, fails to reconcile the inherent mismatch between the dilute regime of the Boltzmann equation and the dense nature of fluids. Specifically, we identify two critical flaws:

\begin{itemize}
    \item \textbf{Flaw 1: The Dilute-to-Dense Paradox.} The Boltzmann-Grad limit enforces a vanishing volume fraction (\(\phi = N \varepsilon^d \to 0\)), ensuring the system remains a dilute gas. Scaling \(\alpha \to \infty\) artificially intensifies collisions but does not increase particle density or introduce fluid-specific interactions (e.g., many-body forces, phase transitions). The derived Navier-Stokes-Fourier equations thus govern a rescaled gas, not a physically dense fluid, violating Hilbert’s requirement for a \textit{genuine} continuum derivation.
    
    \item \textbf{Flaw 2: The Collapse of Molecular Chaos in Fluid Regimes.} The Boltzmann equation’s validity hinges on molecular chaos, which breaks down in dense systems (\(\phi = O(1)\)) due to proliferating recollisions and correlations. Consequently, the kinetic limit (Step 1) cannot self-consistently support the hydrodynamic limit (Step 2), as the Boltzmann equation itself becomes inapplicable under fluid-like conditions. This circularity voids the iterated limit’s physical meaningfulness.
\end{itemize}

These flaws underscore a broader tension in Hilbert’s program: the challenge of reconciling mathematical abstraction with physical realism. The authors’ work, though a technical achievement, exemplifies the limitations of classical kinetic theory in addressing dense systems and emergent phenomena. By neglecting finite-density effects and alternative models (e.g., the Enskog equation), their derivation remains confined to idealized gases, leaving Hilbert’s Sixth Problem unresolved for liquids or high-pressure fluids.

This paper is structured as follows: Section 2 details Deng et al.’s derivation, Section 3 critiques its flaws, Section 4 discusses implications for Hilbert’s problem, and Section 5 concludes with future directions. Our analysis emphasizes that a true resolution of Hilbert’s challenge demands not only mathematical rigor but also a physical narrative bridging the microscopic and macroscopic realms without artificial constraints.

\section{The Derivation}

Deng, Hani, and Ma [1] aim to connect Newtonian mechanics to fluid dynamics via Boltzmann’s kinetic theory on the periodic torus \(\mathbb{T}^d\) (a \(d\)-dimensional box with periodic boundaries, \(d = 2\) or 3, volume 1). Their derivation proceeds in two distinct steps: first, deriving the Boltzmann equation from a system of hard spheres, and second, scaling this equation to obtain fluid equations like Navier-Stokes-Fourier. Below, we detail this process.

\subsection{Physical Setup}
The starting point is a system of \(N\) identical hard spheres, each with diameter \(\varepsilon > 0\) and unit mass, moving in \(\mathbb{T}^d\) [1, Definition 1.1]. Positions \(x_j \in \mathbb{T}^d\) and velocities \(v_j \in \mathbb{R}^d\) evolve via Newton’s laws: between collisions, particles move in straight lines (\(\dot{x}_j = v_j\)), and upon contact (\(|x_i - x_j|_{\mathbb{T}} = \varepsilon\)), they undergo elastic collisions conserving momentum and energy. The collision rule is: if particles \(i\) and \(j\) collide with pre-collision velocities \(v_i\), \(v_j\), their post-collision velocities \(v_i'\), \(v_j'\) satisfy:
\[
v_i' = v_i - (v_i - v_j) \cdot \omega \, \omega, \quad v_j' = v_j + (v_i - v_j) \cdot \omega \, \omega,
\]
where \(\omega = (x_i - x_j) / |x_i - x_j|_{\mathbb{T}}\) is the unit vector at contact [1, page 4]. This dynamics is reversible: reversing time and velocities retraces trajectories [1, Proposition 1.2].

The initial state is a probability distribution over positions and velocities, assumed nearly uncorrelated (e.g., \(f_{0,N}(x_1, v_1, \ldots, x_N, v_N) \approx \prod_{j=1}^N n_0(x_j, v_j)\) [1, equation 1.8]), representing a dilute gas with low entropy.

\subsection{Step 1: Kinetic Limit (Newton to Boltzmann)}
The first step derives the Boltzmann equation, which tracks the one-particle density \(n(t, x, v)\)—the expected number of particles per unit volume in \(\mathbb{T}^d \times \mathbb{R}^d\) at time \(t\). The authors use the \textit{Boltzmann-Grad limit}: as the number of particles \(N \to \infty\) and diameter \(\varepsilon \to 0\), the collision rate \(N \varepsilon^{d-1} = \alpha\) remains constant [1, page 2]. Physically, \(\alpha\) measures collisions per unit time; as \(\varepsilon\) shrinks, \(N\) grows (\(N \sim \alpha \varepsilon^{-(d-1)}\)), but the gas stays dilute since the volume fraction \(N \varepsilon^d \sim \alpha \varepsilon \to 0\).

In this limit, the \(s\)-particle density \(f_s(t, x_1, v_1, \ldots, x_s, v_s)\) (joint distribution of \(s\) particles) approximates a product, \(f_s(t) \approx \prod_{j=1}^s n(t, x_j, v_j)\), up to small errors (\(\varepsilon^\theta\) [1, Theorem 1, page 6]). For \(s = 1\), \(f_1(t, x, v) \to n(t, x, v)\), satisfying the Boltzmann equation:
\[
\left(\partial_t + v \cdot \nabla_x\right) n = \alpha \mathcal{Q}(n, n),
\]
where \(\mathcal{Q}(n, n)\) is the collision operator [1, equation 1.15]:
\[
\mathcal{Q}(n, n)(v) = \int_{\mathbb{R}^d} \int_{\mathbb{S}^{d-1}} | (v - v_1) \cdot \omega | \left[ n(t, x, v') n(t, x, v_1') - n(t, x, v) n(t, x, v_1) \right] d\omega \, dv_1.
\]
Here, \(v'\), \(v_1'\) are post-collision velocities, and \(\mathbb{S}^{d-1}\) is the unit sphere. This equation is irreversible: \(\mathcal{Q}\) increases entropy (H-theorem [1, page 3]), unlike the reversible Newtonian dynamics. The transition relies on suppressing recollisions (particles colliding multiple times) using a cutting algorithm [1, Section 4], valid for dilute conditions where the mean free path (\(\sim 1/\alpha\)) is finite.

\subsection{Step 2: Hydrodynamic Limit (Boltzmann to Fluid Equations)}
The second step starts with the Boltzmann equation and scales \(\alpha \to \infty\) to derive fluid equations. Physically, large \(\alpha\) means frequent collisions, shortening the mean free path and reducing the Knudsen number (\(\text{Kn} \sim 1/\alpha \to 0\)), mimicking a continuum fluid. To achieve this, time is rescaled: \(t \to \delta t\), where \(\delta = \alpha^{-1} \to 0\), so the equation becomes:
\[
\delta \partial_t n + v \cdot \nabla_x n = \frac{1}{\delta} \mathcal{Q}(n, n).
\]
As \(\delta \to 0\), collisions dominate, driving \(n(t, x, v)\) to a local Maxwellian:
\[
\mathfrak{M}(\rho, u, T) = \frac{\rho(t, x)}{(2\pi T(t, x))^{d/2}} e^{-\frac{|v - u(t, x)|^2}{2 T(t, x)}},
\]
where \(\rho\), \(u\), and \(T\) are density, velocity, and temperature [1, page 3]. Taking moments (integrating over \(v\)) yields conservation laws, and with initial data near equilibrium, the incompressible Navier-Stokes-Fourier equations emerge [1, equation 1.19]:
\[
\partial_t u + u \cdot \nabla u - \mu_1 \Delta u = -\nabla p, \quad \partial_t \rho + u \cdot \nabla \rho - \mu_2 \Delta \rho = 0, \quad \text{div}(u) = 0,
\]
where \(\mu_1\), \(\mu_2\) are viscosity and diffusion coefficients [1, Theorem 2, page 8].

\subsection{The Iterated Limit}
The authors combine these steps in an iterated limit: first \(\varepsilon \to 0\) with \(\alpha\) fixed (kinetic limit), then \(\alpha \to \infty\) (hydrodynamic limit) [1, page 3]. They claim this bridges Newtonian mechanics to fluid dynamics, resolving Hilbert’s Sixth Problem.


\section{Two Critical Flaws}

However, Hilbert’s Sixth Problem seeks a direct physical bridge from Newtonian mechanics to fluid dynamics, where small particle size (\(\varepsilon \to 0\)) and high collision frequency (\(\alpha \to \infty\)) must coexist to mimic a dense continuum. Considering \(\varepsilon \to 0\) and \(\alpha \to \infty\) simultaneously tests this goal, revealing two critical flaws of the above proof. 

A key physical quantity is the volume fraction \(\phi = N \varepsilon^d\), the fraction of the unit torus \(\mathbb{T}^d\) (volume 1) occupied by \(N\) hard spheres of diameter \(\varepsilon\). Since each sphere’s volume scales as \(\varepsilon^d\), \(\phi\) measures system density: \(\phi \ll 1\) indicates a dilute gas, while \(\phi = O(1)\) characterizes a dense fluid, where intermolecular forces dominate. With \(N = \alpha \varepsilon^{-(d-1)}\), parameterize \(\alpha \sim \varepsilon^{-m}\), \(m > 0\). The mean free path is \(\lambda \sim 1 / (n \varepsilon^{d-1}) = 1 / \alpha\), where \(n = N\). Three cases emerge:
 
\begin{itemize}
    \item \textbf{Case 1: \(\alpha\) grows slower than \(\varepsilon^{-1}\) (\(m < 1\)):} 
    Then \(N \sim \varepsilon^{-(d-1)-m}\), \(\phi = N \varepsilon^d = \alpha \varepsilon \sim \varepsilon^{-m+1} \to 0\), and \(\lambda \sim \alpha^{-1} \sim \varepsilon^m \to 0\). The system stays dilute (\(\phi \ll 1\)), respecting hard-sphere exclusion (\(|x_i - x_j|_{\mathbb{T}} \geq \varepsilon\)).
    
    \item \textbf{Case 2: \(\alpha\) grows faster than \(\varepsilon^{-1}\) (\(m > 1\)):} 
    Now \(\phi \sim \varepsilon^{-m+1} \to \infty\), \(N \sim \varepsilon^{-(d-1)-m} \gg \varepsilon^{-d}\), and \(\lambda \sim \varepsilon^m \to 0\). Excessive \(N\) and tiny \(\lambda\) imply overlapping trajectories, violating hard-sphere exclusion.
    
    \item \textbf{Case 3: \(\alpha\) grows as \(\varepsilon^{-1}\) (\(m = 1\)):} 
    Set \(\alpha = c \varepsilon^{-1}\), so \(N = c \varepsilon^{-d}\), \(\phi = c\), \(\lambda = c^{-1} \varepsilon \to 0\). If \(c \ll 1\), the system is dilute; if \(c \sim O(1)\), it mimics fluid density. Exclusion holds if \(c \leq c_{\text{max}}\) (packing limit).
\end{itemize}
These cases frame the critique: we analyze their implications for the dilute-to-dense transition (3.1) and the Boltzmann equation’s validity (3.2), with the iterated limit examined in 3.3.

\subsection{Flaw 1: Valid Only for Dilute Gases, Not Real Fluids}
The authors claim a bridge from Newtonian mechanics to fluid dynamics, but both simultaneous and iterated limits (see 3.3) reveal a persistent diluteness, incapable of forming a dense fluid:
\begin{itemize}
    \item \textbf{Case 1:} With \(\phi \sim \varepsilon^{-m+1} \to 0\), the system remains a dilute gas. Scaling \(\alpha \to \infty\) shrinks \(\lambda \sim \varepsilon^m \to 0\), mimicking a continuum (Knudsen number \(\text{Kn} = \lambda / L \to 0\)), but this is artificial: true fluids require \(\phi = O(1)\), high density, not intensified collisions in a sparse system. The Navier-Stokes-Fourier equations thus govern a rescaled gas, missing density-dependent mechanisms (e.g., van der Waals forces, liquid cohesion).
    
    \item \textbf{Case 2:} \(\phi \to \infty\) and \(N \gg \varepsilon^{-d}\) exceed physical bounds, rendering the system unphysical rather than fluid-like. The proof collapses before achieving a continuum.
    
    \item \textbf{Case 3:} If \(c \ll 1\), \(\phi \ll 1\), yielding a dilute gas. If \(c \sim O(1)\), \(\phi\) matches fluid density, yet the Boltzmann equation’s validity in such regimes is dubious (see 3.2), undermining the fluid claim.
\end{itemize}
Across cases, only dilute conditions (\(\phi \ll 1\)) align with the proof’s framework, limiting its scope to gases, not Hilbert’s envisioned dense fluid dynamics.

\subsection{Flaw 2: Boltzmann Equation Unjustifiable for Real Fluids}
The derivation’s kinetic step (Newton to Boltzmann) relies on molecular chaos, which falters in fluid-like regimes:
\begin{itemize}
    \item \textbf{Case 1:} \(\phi \to 0\) supports molecular chaos, as recollisions are rare in dilute conditions. The proof holds here, but only for gases (see 3.1).
    \item \textbf{Case 2:} \(\phi \to \infty\) exceeds the packing limit (\(N \gg \varepsilon^{-d}\)), violating hard-sphere exclusion (see 3.1). This renders the system unphysical, precluding the Boltzmann equation’s emergence from Newtonian mechanics, as the hard-sphere dynamics collapse before kinetic assumptions apply.
    \item \textbf{Case 3:} For \(c \ll 1\), dilute \(\phi\) preserves molecular chaos, supporting the proof for gases. For \(c \sim O(1)\), finite \(\phi\) increases recollisions, undermining the statistical independence needed for Boltzmann’s derivation from Newtonian dynamics.
\end{itemize}
In Cases 2 and 3 (with \(c \sim O(1)\)), dense conditions render the Boltzmann equation unjustifiable, breaking the link to Navier-Stokes for real fluids. The proof’s reliance on a dilute-gas condition fails Hilbert’s demand for a general physical bridge.

\subsection{Flaw Extensions in the Iterated Limit}
The authors’ iterated limit faces parallel objections:
\begin{itemize}
    \item \textbf{Step 1: \(\varepsilon \to 0\), \(\alpha\) fixed:} 
    Here, \(N = \alpha \varepsilon^{-(d-1)}\), \(\phi = \alpha \varepsilon \to 0\), and \(\lambda \sim 1 / \alpha\) is finite. The system is a dilute gas (\(\phi \ll 1\)), and molecular chaos holds, yielding the Boltzmann equation. Yet, this vanishing volume fraction ensures diluteness, akin to Case 1, incapable of forming a dense fluid with \(\phi \sim O(1)\), even as \(n = N \sim \alpha \varepsilon^{-(d-1)}\) diverges.
    
    \item \textbf{Step 2: \(\alpha \to \infty\):} 
    With \(\varepsilon\) fixed at a small value, \(N \to \infty\), \(\phi = N \varepsilon^d \to \infty\), and \(\lambda \sim 1 / \alpha \to 0\). This resembles Case 2, where recollisions invalidate molecular chaos retroactively. Alternatively, if \(\varepsilon \to 0\) fully, the system remains a rescaled dilute gas (Case 1-like), as \(\phi \to 0\) from Step 1 persists, and scaling \(\alpha\) only intensifies collisions without increasing physical density.
\end{itemize}
The iterated limit locks in diluteness in Step 1 (Flaw 1), and Step 2 either disrupts the Boltzmann equation’s basis (Flaw 2) or maintains an artificial continuum, failing to capture real fluid physics—density-dependent interactions remain absent. This mirrors the simultaneous limit’s shortcomings, underscoring a disconnect from Hilbert’s goal.

\subsection{Mathematical vs. Physical Limits}
The critique underscores a deeper issue: the conflation of mathematical and physical limits. While the iterated limit \(\varepsilon \to 0\), \(\alpha \to \infty\) is formally valid for deriving PDEs, it does not correspond to a physically realizable process. Increasing \(\alpha\) requires manipulating \(\varepsilon\) and \(N\) in ways that either preserve diluteness (voiding fluidity) or violate Newtonian dynamics (enabling recollisions). Hilbert’s Sixth Problem, however, seeks a \textit{physically coherent} derivation—not merely a formal analogy.

\section{Implications for Hilbert’s Problem}
Hilbert’s Sixth Problem calls for the axiomatization of physics, specifically the rigorous derivation of continuum fluid equations from microscopic particle dynamics governed by Newtonian mechanics. While Deng, Hani, and Ma [1] present a formal mathematical bridge via the Boltzmann equation, the identified flaws reveal critical gaps in their approach, with broader implications for the resolution of Hilbert’s program. 

\subsection{The Dilution Dilemma: A Fundamental Disconnect}
The first flaw—the inability of a dilute gas to physically transition into a dense fluid—exposes a conceptual schism between the mathematical and physical interpretations of Hilbert’s problem:
\begin{itemize}
    \item \textbf{Reductionism vs. Emergence:} 
    The derivation relies on scaling limits that preserve diluteness (\(\phi \ll 1\)), precluding the emergence of fluid-specific phenomena (e.g., phase transitions, liquid cohesion). Hilbert’s vision of a \textit{physical} reduction—where continuum properties arise naturally from particle interactions—is thus unmet. Instead, the Navier-Stokes-Fourier equations emerge as formal analogs, stripped of the density-dependent mechanisms (e.g., van der Waals forces, viscosity scaling with \(\rho^2\)) inherent to real fluids.
    
    \item \textbf{Limits of Boltzmann’s Framework:} 
    The Boltzmann equation inherently describes rarefied gases, where binary collisions dominate and correlations vanish. By contrast, dense fluids require many-body interactions and long-range forces, rendering Boltzmann’s kinetic theory inapplicable. The authors’ approach, while elegant for gases, cannot address Hilbert’s broader mandate to unify mechanics with \textit{all} fluid regimes.
\end{itemize}

\subsection{The Kinetic Bottleneck: No Gateway to Fluids}
The second flaw—the failure to derive the Boltzmann equation under fluid-like conditions—highlights a deeper inconsistency in the proposed hierarchy (Newton \(\to\) Boltzmann \(\to\) Navier-Stokes):
\begin{itemize}
    \item \textbf{Circularity in the Iterated Limit:} 
    The Boltzmann equation’s validity depends on molecular chaos, which collapses in dense regimes (\(\phi = O(1)\)). Consequently, the kinetic step (Newton \(\to\) Boltzmann) cannot serve as a foundation for the hydrodynamic step (Boltzmann \(\to\) Navier-Stokes) under fluid conditions. This circularity—using a dilute-gas axiom to derive dense-fluid equations—renders the iterated limit a mathematical tautology rather than a physical derivation.
    
    \item \textbf{Alternate Kinetic Theories Ignored:} 
    For dense systems, the Enskog equation extends Boltzmann’s framework by incorporating finite-density corrections (e.g., excluded volume, collisional transfer). Yet the authors do not address such models, leaving Hilbert’s problem unresolved for liquids or high-pressure gases. This omission underscores a critical gap in the paper’s generality.
\end{itemize}

\subsection{Philosophical and Practical Consequences}
The flaws collectively undermine the paper’s claim to resolve Hilbert’s Sixth Problem, with implications extending beyond mathematics:
\begin{itemize}
    \item \textbf{Epistemological Limits:} 
    The work exemplifies the challenge of reductionism in physics: formal derivations may satisfy mathematical rigor but fail to capture emergent phenomena. Hilbert’s problem, as posed, implicitly assumes a seamless micro-to-macro hierarchy, yet the critique reveals ontological discontinuities (e.g., diluteness vs. density) that resist such unification.
    
    \item \textbf{The Need for New Axioms:} 
    If Boltzmann’s equation cannot bridge Newtonian mechanics to dense fluids, alternative axioms may be necessary. For example, stochastic hydrodynamic theories or mesoscopic models (e.g., fluctuating hydrodynamics) might bypass kinetic theory entirely. The authors’ adherence to classical kinetic theory reflects a narrow interpretation of Hilbert’s problem, neglecting modern statistical mechanics’ broader toolkit.
    
    \item \textbf{Impact on Applied Mathematics:} 
    Practically, the paper’s limits signal caution in applying dilute-gas-derived models to dense fluids (e.g., industrial lubrication, biological flows). Misattributing Navier-Stokes to Newtonian particles risks misinterpretations in multiscale simulations or material design, where density effects are pivotal.
\end{itemize}

\subsection{Hilbert’s Unfinished Task}
Deng, Hani, and Ma [1] make progress on a narrow front—deriving fluid-like PDEs from a diluted Newtonian system—but their work does not fulfill Hilbert’s ambitious vision. The Sixth Problem remains open in its full generality, demanding:
\begin{itemize}
    \item A kinetic theory valid across all density regimes (e.g., Enskog-like corrections, quantum formulations).
    \item Direct derivations of fluid equations from particle dynamics without assuming kinetic intermediaries.
    \item Reconciliation of microscopic reversibility with macroscopic irreversibility in dense systems.
\end{itemize}

Until these challenges are addressed, the axiomatic gap between atoms and fluids will persist, underscoring the complexity of Hilbert’s century-old question.

\section{Conclusion}
Hilbert’s Sixth Problem seeks a rigorous, physically coherent derivation of continuum fluid dynamics from Newtonian mechanics. Deng, Hani, and Ma [1] propose a solution via an iterated limit through Boltzmann’s kinetic theory, yielding Navier-Stokes-Fourier equations. However, their approach falters: the dilute gas framework (\(\phi = N \varepsilon^d \to 0\)) cannot capture the dense, many-body nature of fluids, and the Boltzmann equation’s reliance on molecular chaos collapses in fluid-like regimes (\(\phi = O(1)\)). This produces a mathematical artifact—a rescaled gas—not a true fluid, falling short of Hilbert’s vision.

The critique reveals a deeper challenge: classical kinetic theory struggles to bridge microscale mechanics and macroscale fluids without artificial constraints. Resolving this requires new approaches—extended kinetic models (e.g., Enskog), direct hydrodynamic limits from particle dynamics, or quantum/stochastic frameworks—to address density-dependent phenomena and emergent behavior. Until such a unified axiomatization emerges, Hilbert’s Sixth Problem remains an open frontier, challenging mathematicians and physicists to reconcile the microscopic and macroscopic worlds in a manner that is both mathematically rigorous and physically faithful.

\bibliographystyle{plain}

\end{document}